\documentclass[a4paper,11pt]{article}

\usepackage[dvipsnames]{xcolor}
\usepackage[english]{babel}
\usepackage{amssymb, enumerate}
\usepackage{latexsym}
\usepackage[all,cmtip]{xy}
\usepackage{amsmath,amssymb,amsthm}
\usepackage{comment}
\usepackage{tikz}
\usepackage{marvosym}
\usepackage{tabu}

\usepackage{color}

\usepackage{tikz}
\usepackage{tikz-cd}
\usepackage{tkz-euclide}
\usepackage{caption}
\usepackage{enumerate}
\usepackage[all]{xy}
\usepackage{MnSymbol}
\usepackage[margin=1.59in]{geometry}
\usepackage{float}
\usepackage{thmtools}
\usepackage{multicol}
\restylefloat{table}
\usepackage{footmisc}
\renewcommand{\thefootnote}{\fnsymbol{footnote}}
\usetikzlibrary{shapes.geometric,decorations}
\usepackage{titlesec}
\usepackage{listings}
\usepackage{relsize}
\definecolor{burgundy}{rgb}{0.5, 0.0, 0.13}
\definecolor{asparagus}{rgb}{0.163, 0.66, 0.42}
\usepackage[colorlinks=true, citecolor=Blue, linkcolor=Blue, urlcolor=Blue, linktocpage=true, pdfpagelabels, ]{hyperref}
\usepackage{cleveref}

\newtheorem{theorem}{Theorem}[section]
\newtheorem*{theorem*}{Theorem}

\newtheorem{lemma}[theorem]{Lemma}
\newtheorem{corollary}[theorem]{Corollary}
\newtheorem{proposition}[theorem]{Proposition}
\newtheorem{remark}[theorem]{Remark}

\newcommand{\stab}{\operatorname{Stab}}
\newcommand{\aut}{\operatorname{Aut}}
\newcommand{\out}{\operatorname{Out}}

\newcommand{\PI}{\mathbb{Z}\operatorname{Irr}}
\newcommand\blfootnote[1]{%
	\begingroup
	\renewcommand\thefootnote{}\footnote{#1}%
	\addtocounter{footnote}{-1}%
	\endgroup
}
\begin{document}\setcounter{MaxMatrixCols}{50}
\title{$2$-Blocks whose defect group is homocyclic and whose inertial quotient contains a Singer cycle II}
	\author{Elliot Mckernon$^{\dagger}$}
	\date{\today}
	\maketitle
\begin{abstract} We consider $2$-blocks of finite groups with defect group $D=Q \times R$ and inertial quotient $\mathbb{E}$ where $Q \cong (C_{2^m})^n$, $R \cong C_{2^r}$, and $\mathbb{E}$ contains a Singer cycle of $\aut(Q)$ (an element of order $2^n-1$). We classify such blocks up to Morita equivalence when either $\mathbb{E}$ is cyclic or $r=1$. We achieve a partial classification when $r>1$ and $E$ is non-cyclic. 
		\blfootnote{${}^{\dagger}$ School of Mathematics, University of Manchester, Manchester, M13 9PL, United Kingdom.\\ Email:
			elliot.mckernon@manchester.ac.uk}

		
	\end{abstract}
\section{Introduction}

Let $p$ be a prime number and $\mathcal{O}$ a complete discrete valuation ring with algebraically closed residue field $k$ and field of fractions $K$ of characteristic zero. We assume that $K$ contains a $|G|$-th primitive root of unity for any finite group $G$ we consider.

In a previous paper, \cite{Elliot}, we classified $2$-blocks with defect group $(C_{2^m})^n$ whose inertial quotient contained a Singer cycle. In this paper we consider similar blocks with the addition of a cyclic factor to the defect group, centralised by the inertial quotient. We will prove \Cref{Main2} and \Cref{Maink}. 

Recall that a block of $G$ with defect group $D$ is said to be \textit{inertial} if it is basic Morita equivalent to its Brauer correspondent in $N_G(D)$.

 \begin{theorem} \label{Main2} Let $G$ be a finite group and $B$ a block of $\mathcal{O}G$ with defect group $D=Q \times R$ and inertial quotient $\mathbb{E}$ such that $Q \cong (C_{2^m})^n$, $R \cong C_{2^r}$, and $\mathbb{E}$ centralises $R$ and contains an element of order $2^n-1$, with $n,m,r \geq 1$. Then $\mathbb{E}=E \rtimes F$, where $E \cong C_{2^n-1}$ and $F$ is cyclic of order dividing $n$, and one of the following occurs:
 	
 	\begin{enumerate} 
 \item $B$ is Morita equivalent to the principal block of one of: 	\begin{enumerate}
 	\item $(Q \rtimes (E \rtimes F)) \times R$. \label{2Inert}
 	\item $(SL_2(2^n) \rtimes F) \times R$ (only if $m=1$). \label{2SL}
 	\item $J_1 \times R$ (only if $m=1$ and $n=3$). \label{2J1}
 \end{enumerate} 
\item \label{2coverinert} $B$ covers an inertial block of some $N \lhd G$ satisfying $[G:N]=[R:R \cap N]$, and if $\mathbb{E}$ is cyclic, then $B$ lies in case (\ref{2Inert}).
 	\end{enumerate} 
 \end{theorem}
 
\begin{theorem}	\label{Maink}
	Let $\tilde{B}$ be a block of $kG$ with defect group $D=Q \times R$ and inertial quotient $\mathbb{E}$ as in the hypothesis of \Cref{Main2}. Further suppose that $R \cap E(G)=1$ or $R \leq E(G)$ (which is implied if $r=1$). Then $\tilde{B}$ is Morita equivalent to the principal block over $k$ of one of (\ref{2Inert}), (\ref{2SL}), or (\ref{2J1}) from \Cref{Main2}.  \end{theorem}

To proceed, we will collect standard Clifford theoretic results and other tools in \Cref{section tools}. In \Cref{section IQ}, we consider the relationship between the inertial quotients of a block covering another when the normal subgroup has $p$-power index. In \Cref{section summary} we summarise the salient results from \cite{Elliot}, before proving \Cref{Main2} and \Cref{Maink} in \Cref{section proof}.

\section{Collecting Tools} \label{section tools}
Suppose that $G$ is a finite group and $N \lhd G$. We say a block algebra $B$ of $\mathcal{O}G$ \textit{covers} a block $b$ of $\mathcal{O}N$ if $Bb \neq 0$. This relates the structures of $B$ and $b$:

\begin{proposition} (\cite[15.1]{alperin}) \label{alp15} Let $N \lhd G$ be finite groups, and let $B$ be a $p$-block of $\mathcal{O}G$ covering a block $b$ of $\mathcal{O}N$. Then the following hold: \begin{enumerate} \item The blocks of $N$ covered by $B$ form a $G$-conjugacy class. 
		\item Each defect group of $b$ is the intersection of a defect group of $B$ with $N$. \label{alpintersect}
		\item There is a block $B'$ of $G$ covering $b$ such that $B'$ has a defect group $D'$ satisfying 
		
		{\centering $ [D':D' \cap N]=[\operatorname{Stab}_G(b):N]_p.$ \label{alpblockabove}\par} 
		
		\item If the centraliser in $G$ of a defect group of is contained in $N$ then $B=b^G$ and $B$ is the only block of $G$ covering $b$. \label{alp CG(D)<N}	\end{enumerate}\end{proposition}
\color{black}
Comparing $B$ and $b$ is particularly effective when we have additional information about $G/N$, as we show in the next few lemmas. 

\begin{lemma}  \label{index p}
Let $G$ be a finite group and $N \lhd G$ such that $[G:N]$ is a power of $p$. Let $B$ be a $p$-block of $\mathcal{O}G$ covering a block $b$ of $\mathcal{O}N$. Then $B$ is the unique block of $G$ covering $b$. If $b$ is $G$-stable, then $B$ and $b$ share a block idempotent and $G=DN$. Further, $b$ is $G$-stable if and only if $[D:D \cap N]=[G:N]$. \end{lemma}
\begin{proof}
	That $B$ is unique is \cite[V., Lemma 3.5.]{feit}, and \cite[Proposition 6.8.11]{linckelmann} tells us that $B$ and $b$ share a block idempotent. \Cref{alp15} (\ref{alpintersect}) tells us that $[D:D \cap N]=[\operatorname{Stab}_G(b):N]_p$. Since $[G:N]$ is a power of $p$, $|DN/N|=|\operatorname{Stab}_G(b)/N|$, implying $|DN|=|\operatorname{Stab}_G(b)|$. Thus, if $b$ is $G$-stable, we have $[D:D \cap N]=[G:N]$, and $G=DN$. If we have $[D: D \cap N]=[G:N]$ then $[\stab_G(b):N]=[G:N]$, implying $\stab_G(b)=G$ and that $b$ is $G$-stable.  
\end{proof}

\begin{lemma} (\cite[Lemma 2.4]{CesC2^5}) \label{solv quotient}
	Let $N \lhd G$ be finite groups such that $G/N$ is solvable, and $B$ a quasiprimitive block of $G$ with abelian defect group $D$. Then $DN/N$ is a Sylow-$p$-subgroup of $G/N$. 
\end{lemma}

If $A$ and $B$ are $\mathcal{O}$-algebras, they are Morita equivalent if their module categories are equivalent. This occurs if and only if there is an $(A,B)$-bimodule $M$ and a $(B,A)$-bimodule $N$ such that $M \otimes_B N \cong A$ as $(A,A)$-bimodules and $N \otimes_A B \cong B$ as $(B,B)$-bimodules. If $N$ and $M$ have endopermutation source, $A$ and $B$ are said to be \textit{basic Morita equivalent}. If they have trivial source, they are \textit{Puig equivalent}.

A $B$-subpair is a pair $(P,B_P)$ where $P$ is a $p$-subgroup of $G$ and $B_P$ is a block of $\mathcal{O}PC_G(P)$ with Brauer correspondent $B$. When $D$ is a defect group of $B$, the $B$-subpairs $(D,B_D)$ are $G$-conjugate, and we write $N_G(D,B_D)$ for the stabiliser in $N_G(D)$ of $B_D$. The \textit{inertial quotient} of $B$ is then $E=N_G(D,B_D)/DC_G(D)$, a $p'$-group unique up to isomorphism. If $A$ and $B$ are basic Morita equivalent blocks, then they have the same defect group and inertial quotient \cite{linckelmann}.

When $B$ has abelian defect group, $B$ is \textit{nilpotent} if and only if its inertial quotient is trivial, and \textit{nilpotent-covered} if there is a nilpotent block covering $B$. 

\begin{proposition} \label{zhou}
Let $G$ be a finite group and $N \lhd G$, and let $B$ be a block of $\mathcal{O}G$ covering a block $b$ of $\mathcal{O}N$. \begin{enumerate}
		\item If $b$ is nilpotent-covered, then $b$ is inertial. \label{zhou b nilp b inert}
		\item If $p \not \divides [G:N]$ and $b$ is inertial, then $B$ is inertial.  \label{zhou b inert B inert}
	\end{enumerate}
\end{proposition} 
\begin{proof}
	(\ref{zhou b nilp b inert}) is \cite[Corollary 4.3]{Puigext2}, and (\ref{zhou b inert B inert}) is the main result of \cite{zhou}. 
\end{proof}

Comparing blocks when $[G:N]$ is a power of $p$ is complex, but we do have the following result due to Koshitani and K\"ulshammer:

\begin{theorem} \label{koshkul} (\cite{koshitanikulshammer})
Let $G$ be a finite group and $N \lhd G$ such that $G/N$ is a $p$-group. Let $b$ be a $G$-stable $p$-block of $kN$, and let $B$ the unique block of $kG$ covering $b$. If $B$ has an abelian defect group $D$ such that $D \cap N$ has a complement $R$ in $D$, then $B \cong kR \otimes_{k} b$ as $k$-algebras. 
\end{theorem}

In \cite{watanabesplittingtheorem}, Watanabe extends this to $\mathcal{O}$ when a certain perfect isometry can be found. The $*$ in the following denotes the $*$-construction of characters due to Brou\'e and Puig \cite{star}. 

\begin{proposition} (\cite[Lemma 3]{watanabesplittingtheorem})\label{watanabesplitting}
	Let $N  \lhd G$ be finite groups and $B$ a block of $\mathcal{O}G$ covering a $G$-stable block $b$ of $\mathcal{O}N$. Suppose $B$ has abelian defect group $D$ such that $D=(D \cap N) \times R$ with $G=N \rtimes R$. Let $B_D$ be a Brauer correspondent of $B$ in $C_G(D)$, and set $B_R:=(B_D)^{C_G(R)}$. If there's a perfect isometry $\epsilon:\PI(C_G(R),B_R) \to \PI(G,B)$ such that $\epsilon(\lambda * \zeta)=\lambda * \epsilon(\zeta)$ for every $\lambda \in \operatorname{Irr}(R)$ and $\zeta \in \operatorname{Irr}(B_R)$, then $B \cong b \otimes_{\mathcal{O}} \mathcal{O}R$ as $\mathcal{O}$-algebras.
\end{proposition}

 There's such no tool over $\mathcal{O}$ or $k$ when $G$ is a non-split extension, though this can be managed in certain circumstances (e.g. in \cite[\S 2]{eatonlivesrank3}). When the $E$ is cyclic and $B$ has $|E|$ simple modules, we can apply another result of Watanabe to find the perfect isometry required by \Cref{watanabesplitting}.  Recall that, for a block $B$ of a group $G$, $l(B)$ denotes the number of irreducible characters of $G$ that lie in $B$. 

\begin{theorem} (\cite{watanabepis}) \label{watanabepis} 
	Let $G$ be a finite group and $B$ a block of $G$ with defect group $D$ and cyclic inertial quotient $E$. Let $(D,B_D)$ be a maximal $(G,B)$-Brauer pair, and set $N:=N_G(D,B_D)$. If $l(B)=|E|$, and $C_E(d)=1$ for each $d \in [N,D] \backslash \{1\}$, then there exists a perfect isometry $\iota:\PI(N,{B_D}^N) \to \PI(G,B)$ such that $\iota(\lambda * \zeta)=\lambda * \iota(\zeta)$ for $\lambda \in \operatorname{Irr}(C_D(N))$ and $\zeta \in \PI(N,{B_D}^N)$. 
\end{theorem}

%
%

When the defect group of a block is normal in $G$, we understand the block structure well using the following result of K\"ulshammer:

\begin{lemma} \label{kuls} (\cite{kulscrossed})
	Let $N \lhd G$ be finite groups, and let $B$ be a block of $\mathcal{O}G$ with defect group $D$ and inertial quotient $E$. If $D \lhd G$, then $B$ is Puig equivalent to a twisted group algebra $\mathcal{O}_\gamma(D \rtimes E)$, where $\gamma \in O_{p'}(H^2(E,\mathcal{O}^\times))$.
\end{lemma}

 This will be particularly useful for us, since the inertial quotients under consideration have trivial Schur multipliers by \cite[Proposition 5.3]{Elliot}, and so when $D \lhd G$ we will have that $B$ is Puig equivalent to $\mathcal{O}(D \rtimes E)$. 

When proving classifications up to Morita equivalence, there are a number of useful reductions we can apply, which give us blocks with a simpler structure that lie in the same Morita equivalence class. We start with the first Fong reduction:
\begin{theorem} (\cite[6.8.3]{linckelmann}) \label{F1} 
Let $G$ be a finite group with $N \lhd G$, and let $B$ a block of $\mathcal{O}G$ covering a block $b$ of $\mathcal{O}N$. Then there is a unique block $\tilde{B}$ of $\operatorname{Stab}_G(b)$ that is covered by $B$ and covers $b$. Further, $B$ and $\tilde{B}$ have the same defect group and fusion system.
\end{theorem}

As described in \cite[Proposition 2.2]{eaton3}, one can derive the following from \cite{Puigext}:

\begin{theorem} \label{Fong in 24} 

	Let $G$ be a finite group with $N \lhd G$, and $B$ a block of $G$ with defect group $D$ covering a $G$-stable nilpotent block of $N$ with defect group $D \cap N$. Then there is a finite group $L$ and $M \lhd L$ such that: $M \cong D \cap N$ and $L/M \cong G/N$; there is a subgroup $D_L \leq L$ with $D_L \cong D$ and $D_L \cap M \cong D \cap N$; and there is a central extension $\tilde{L}$ of $L$ by a $p'$-group, and a block $\tilde{B}$ of $\mathcal{O}\tilde{L}$ which is Morita equivalent to $B$ and has defect group $\tilde{D} \cong D_L \cong D$. 
\end{theorem}

Note that \cite[6.8.13]{linckelmann} tells us the Morita equivalences in Theorems \ref{F1} and \ref{Fong in 24} are basic.

\begin{remark} \label{reduce remark} \normalfont
Let $G$ be a finite group and $B$ a block of $G$. We say $B$ is \textit{quasiprimitive} if, for every $N \lhd G$, $B$ covers a $G$-stable block of $N$. We say $B$ is \text{reduced} if it is quasiprimitive and if $B$ covering a nilpotent block of some $N \lhd G$ implies $N \leq O_2(G)Z(G)$. Applying \Cref{F1} and \Cref{Fong in 24} lets us replace a given block of $G$ with a reduced block of some $\tilde{G}$, which is basic Morita equivalent to the original. Further, this $\tilde{G}$ is smaller than $G$ in a certain sense:
\end{remark}

\begin{corollary} (\cite[Corollary 2.3]{eaton4})\label{F2 cor}
	Let $N \lhd G$ be finite groups such that $N \not \leq O_2(G)Z(G)$. Let $B$ be a quasiprimitive block of $\mathcal{O}G$ covering a nilpotent block $b$ of $N$. Then there's a finite group $H$ with $[H:O_{p'}(Z(H))] < [G:O_{p'}(Z(G))]$ and a block $B_H$ that is Morita equivalent to $B$ and has isomorphic defect group.
\end{corollary}

\section{Preservation of the inertial quotient} \label{section IQ}

Let $G$ be a finite group with $N \lhd G$, such that $[G:N]$ is a power of $p$, and let $B$ be a block of $\mathcal{O}G$ with abelian defect group $D$ and inertial quotient $E$ covering a block of $\mathcal{O}N$. In this section, we prove that if $[D,E] \leq N$, then the inertial quotient of $b$ is isomorphic to $E$. 

We start with a lemma showing that conjugates of Brauer corresponding blocks are the Brauer correspondents of conjugates of those blocks.

\begin{lemma} \label{brauer corre conj}
	Let $H \leq G \leq X$ be finite groups, and $B$ a block of $G$ with defect group $D$ such that $C_G(D) \leq H$. If $b$ is a block of $H$ satisfying $b^G=B$, then $({}^xb)^G={}^xB$, for all $x \in X$. In particular, $\operatorname{Stab}_X(b) \leq \operatorname{Stab}_X(B)$. 
\end{lemma}
\begin{proof} 
	Let $e_B$ and $e_b$ be the block idempotents of $B$ and $b$ in $\mathcal{O}G$ and $\mathcal{O}H$ respectively, and write $e_B=\sum \limits_{g \in G} \beta_g g$. Recall that $b$ is the Brauer correspondent of $B$ if $e_b\operatorname{Br}_D^G(e_B) \neq 0$. Then $$\operatorname{Br}_D^G({}^xe_B)=\sum \limits_{g \in C_G(D)} \beta_g ({}^xg)={}^x\left( \sum \limits_{g \in C_G(D)} \beta_g g\right) = {}^x\bigl(\operatorname{Br}_D^G(e_B)\bigr).$$ Therefore, $({}^xe_b)\operatorname{Br}_D^G({}^xe_B)=
	({}^xe_b)\left({}^x\operatorname{Br}_D^G(e_B)\right)=
	{}^x\left(e_b\operatorname{Br}_D^G(e_B)\right) \neq 0$, and so $({}^xb)^G={}^xB$. In particular, if ${}^xb=b$, then $({}^xb)^G={}^xB=B$. 
\end{proof}

	Let $G$ be a finite group and $N$ a normal subgroup with index $p^k$ for a prime $p$. Let $b$ be a $p$-block of $N$ covered by a block $B$ of $G$ with abelian defect group $D$ such that $[D,E] \leq D \cap N$ and $[G:N]=[D:D \cap N]$. Then $b$ has inertial quotient isomorphic to $E$.

\begin{theorem} \label{IQ con}
Let $G$ be a finite group with $N \lhd G$ such that $[G:N]$ is a power of $p$. Let $B$ be a block of $\mathcal{O}G$ with abelian defect group $D$ and inertial quotient $E$, and let $b$ be a $G$-stable block of $\mathcal{O}N$ covered by $B$. If $[D,E] \leq N$, then the inertial quotient of $b$ is isomorphic to $E$ and its action on $D$ is identical to that of $E$.
\end{theorem}

Note that the assumption that $[D,E] \leq N$ is redundant if $D \cap N$ has a complement $S$ in $D$ by applying an argument from \cite{koshitanikulshammer}: since $D$ is a $p$-group while $E$ is a $p'$-group, if $D \cap N$ is an $N_G(D,B_D)$-stable direct factor of $D$, then Maschke's theorem (as in \cite[Theorem 3.3.2]{gorenstein}) tells us that $D \cap N$ has an $N_G(D,B_D)$-stable complement in $D$, which we can assume is $S$ in our notation. Further, $[S,N_G(D,B_D)] \leq S$, but since $G/N \cong S$ is abelian, $[G,G] \leq N$, implying that $[S,N_G(D,B_D)] \leq S \cap N =1$. Thus $S \leq C_D(E)$. 

Recall that a subgroup $H \leq G$ is said to \textit{control fusion of $B$-subpairs} if, for each $B$-subpair $(P,B_P)$, $N_G(P,B_P)=N_H(P,B_P)C_G(P)$ \cite{alpbrou}.

\begin{proof}
Since $[G:N]$ is a power of $p$ and $b$ is $G$-stable, \Cref{index p} tells us that $B$ is the unique block covering $b$, that $G=ND$, and that $[G:N]=[D:D \cap N]$. Since $S \leq C_D(E)$, the same relationship holds between blocks of $N_G(D)$ covering blocks of $N_N(D)$, and likewise for $C_G(D)$ and $C_N(D)$. 
	
Let $B_D$ be a Brauer correspondent of $B$ in $C_G(D)$ and write $P:=D \cap N$. It is easy to see that $C_G(D) \leq N_G(D,B_D) \cap C_G(P)$, and, since $N_G(D,B_D)$ centralises $S$, we have $N_G(D,B_D) \cap C_G(P) \leq C_G(S) \cap C_G(P) = C_G(D)$. Thus we have $C_G(D)=N_G(D,B_D) \cap C_G(P)$, and so we can write $$\frac{N_G(D,B_D)}{C_G(D)}=\frac{N_G(D,B_D)}{N_G(D,B_D) \cap C_G(P)}$$ 
Since $P$ is $N_G(D,B_D)$-stable, we also have $N_G(D,B_D) \leq N_G(P)$, implying that $N_G(D,B_D) \leq N_G(C_G(P))$. This lets us apply the second isomorphism theorem, giving us $$\frac{N_G(D,B_D)}{N_G(D,B_D) \cap C_G(P)} \cong \frac{N_G(D,B_D)C_G(P)}{C_G(P)}.$$ 
Next, we claim that $N_G(D,B_D)$ stabilises the respective Brauer correspondents of $b$ in $C_N(P)$ and of $B$ in $C_G(P)$. Note that $C_{C_G(P)}(D)=C_G(D)$, so we can define $B_P:={B_D}^{C_G(P)}$. Next, we choose a block $b_P$ of $C_N(P)$ covered by $B_P$ and satisfying $(b_P)^N=B$. By the first paragraph, $b_P$ and $B_P$ share a block idempotent; $b_P$ is $C_G(P)$-stable; and $B_P$ is the unique block of $C_G(P)$ covering $b_P$. Next, let $x \in N_G(D)$ and consider the block ${}^xB_D$ of $C_G(D)$ and the corresponding block $B'_P:=({}^xB_D)^{C_G(P)}$ of $C_G(P)$. Then $B'_P$ uniquely covers a $C_G(P)$-stable block $b'_P$ of $C_N(P)$ such that $({b'_P})^N=b$. In this situation, \Cref{brauer corre conj} tells us that $B'_P=({}^xB_D)^{C_G(P)}={}^x \bigl( {B_D}^{C_G(P)} \bigr)={}^xB_P$, also implying that $b'_P={}^xb_P$. Thus, if $x \in N_G(D,B_D)$, then we have $B'_P={}^x B_P=B_P$, and so $b'_P={}^xb_P=b_P$. Therefore we have that $N_G(D,B_D) \leq N_G(P,B_P)$.

Next, since $D$ is abelian, Alperin's fusion theorem tells us that $N_G(D,B_D)$ controls fusion of $B$-subpairs. Since $(P,B_P) \leq (D,B_D)$, this means that $N_G(P,B_P)=N_{N_G(D,B_D)}(P,B_P) C_G(P)=$ $\bigl(N_G(D,B_D) \cap N_G(P,B_P) \bigr) C_G(P)$. Since we proved that $N_G(D,B_D) \leq N_G(P,B_P)$, this implies that $N_G(P,B_P)=N_G(D,B_D) C_G(P)$. Therefore we can write $$\frac{N_G(D,B_D)C_G(P)}{C_G(P)} = \frac{N_G(P,B_P)}{C_G(P)}.$$

$B_P$ is the unique block of $C_G(P)$ covering $b_P$, and $b_P$ is $C_G(P)$-stable, telling us that $\operatorname{Stab}_G(B_P)=\operatorname{Stab}_G(b_P)$. Thus, to conclude, we define a map $\varphi:N_N(P,b_P) \to N_G(P,B_P)/C_G(P)$ given by the inclusion into $N_G(P,B_P)$ composed with the canonical projection $g \mapsto gC_G(P)$. By \Cref{index p}, if $g \in N_G(P,B_P)$, then we can write $g=nd$ for $n \in N_N(P,b_P)$ and $d \in D$, so that $\varphi(g)=gC_G(P)=ndC_G(P)=nC_G(P)=\varphi(n)$. Since $n \in N_N(P,b_P)$, $\varphi$ is surjective. Next, observe that $\ker (\varphi)=\{n \in N_N(P,b_P): nC_G(P)=C_G(P) \}=N_N(P,b_P) \cap C_G(P) = C_N(P)$. Thus, we finish the proof by applying the first isomorphism theorem: $$\frac{N_G(P,B_P)}{C_G(P)} \cong \frac{N_N(P,b_P)}{C_N(P)}. $$  \qedhere \end{proof}

Note that \Cref{IQ con} can be deduced from \cite[Lemma 3.6]{pzhou}, though our proof avoids the use of pointed defect groups.

\section{Singer Cycles and Results from \cite{Elliot}} \label{section summary}
\cite[Theorem 1.1]{Elliot} gives the basic Morita equivalence classes of $2$-blocks with defect group $(C_{2^m})^n$ and inertial quotient $\mathbb{E}$ containing an element of order $2^n-1$. We collate several lemmas from that paper in \Cref{Lemma}, but we start by recalling the definition of a Singer cycle.

Let $p$ be a prime, and $P \cong (C_p)^n$. Then $\aut(P) \cong GL_n(p)$ contains elements of order $p^n-1$. These elements are called \textit{Singer cycles}. Further, if $Q \cong (C_{p^m})^n$, then $\aut(Q)/O_p(\aut(Q)) \cong GL_n(p)$, so elements of order $p^n-1$ in $\aut(Q)$ act transitively on the non-trivial elements of $\Omega(Q):=\{q \in Q: q^p=1\} \cong P$, and freely on the non-trivial elements of $Q$.

\begin{proposition} \label{Singers}
Let $n$ be a natural number and $p$ be a prime. Let $Q \cong (C_{p^m})^n$ and define $P:=\Omega(Q) \cong (C_p)^n$ and $G:=\aut(Q)/O_p(\aut(Q))$. Then $G \cong GL_n(p)$ contains a Singer cycle and the following hold:\begin{enumerate}
\item \label{Singers trans}The Singer cycles in $G$ are the elements of maximal order, and the subgroups of $G$ that are generated by a Singer cycle are conjugate in $G$, and act regularly on $P \backslash \{1\}$ and freely on $Q \backslash \{1\}$. 
\item Let $\mathbb{E}$ be a $p'$-subgroup of $G$ containing a Singer cycle. Then either $\mathbb{E}\cong C_{2^n-1}$, or $\mathbb{E} \cong C_{2^n-1} \rtimes F$ where $F$ is cyclic of order dividing $n$. \label{Singers autos}
\item If $H$ is a non-trivial normal subgroup of $\mathbb{E}$, then $C_Q(H)=1$. \label{Singers auto nofix}
\end{enumerate}  
\end{proposition}
\begin{proof}
This follows from \cite[\S 2]{Elliot}. \qedhere	

\end{proof}

In light of \Cref{Singers}, when we consider a block with defect group $(C_{2^m})^n \times C_{2^r}$ whose inertial quotient contains an element of $2^n-1$ and centralises $C_{2^r}$, we understand the action on the defect group well. In \Cref{Lemma} we summarise the consequences of this action, derived from \cite[\S 5]{Elliot}.

Recall the following definitions and facts from \cite[\S 11]{aschbacher}: a group is \textit{quasisimple} if it is a perfect central extension of a simple group. If $G$ is a finite group, a \textit{component} of $G$ is a subnormal, quasisimple subgroup, and the \textit{layer} of $G$, denoted $E(G)$, is the central product of the components of $G$. The \textit{Fitting subgroup} of $G$, $F(G)$, is the product of the $p$-cores of $G$ for each $p$ dividing $|G|$. The \textit{generalised Fitting subgroup}, denoted $F^*(G)$, is the central product of $E(G)$ and $F(G)$. Note that $E(G), F(G)$ and $F^*(G)$ are normal in $G$, and that $C_G(F^*(G)) \leq F^*(G)$, providing an injective homomorphism $G/F^*(G) \to \out(F^*(G))$.

\begin{lemma} \label{Lemma}
Let $G$ be a finite group, and $B$ a reduced $2$-block of $G$ with defect group $D$ and inertial quotient $\mathbb{E}$ such that $D=Q \times R$ with $Q \cong (C_{2^m})^n$ and $R \cong C_{2^r}$, and such that $\mathbb{E}$ centralises $R$ and contains an element of order $2^n-1$. Then the following hold: 
\begin{enumerate}
\item If $N \lhd G$, then $Q \cap N \cong (C_{2^{m'}})^n$ for $m' \leq m$, and $R \cap N \cong C_{2^{r'}}$ for $r' \leq r$. \label{DcapN}
\item If $D$ is not normal in $G$, and $H<D$ such that $H \lhd G$, then $H$ is centralised by a non-trivial subgroup of $\mathbb{E}$. In particular, $O_2(G) \leq R$. \label{O2}
\item If $D$ is not normal in $G$, then $O_{2'}(G) \leq Z(G) \leq Z(F^*(G)) \leq O_2(G)Z(G)$. \label{O2'}
\item If $D$ is not normal in $G$, then the components of $G$ are permuted transitively by $\mathbb{E}$, and $C_D(\mathbb{E}) \cap E(G)$ is contained in the intersection of all the components. \label{comptrans}
\item If $D$ is not normal in $G$, then $\Omega(Q) \leq E(G)$ and there is an embedding $\mathbb{E} \hookrightarrow \operatorname{Aut}(L) \wr S_t$, where $L$ is a component of $G$ and $t$ is the number of components of $G$.  \label{E in out}
\item If $D$ is not normal in $G$, then $t=1$ and $E(G)=L$ is quasisimple. \label{onecomp solv}
\end{enumerate}
\end{lemma}
\begin{proof}
 For (\ref{DcapN}), the structure of $Q \cap N$ follows from \cite[Lemma 5.1]{Elliot}. (\ref{O2}) is \cite[Proposition 4.4]{Elliot}. (\ref{O2'}) follows from \cite[Lemma 4.3]{Elliot} and (\ref{O2}). (\ref{comptrans}) is \cite[Lemma 5.4]{Elliot}. Applying (\ref{DcapN}), either $Q \cap E(G)=1$ or $\Omega(Q) \leq E(G)$. If $Q \cap E(G)=1$, then $B$ covers a nilpotent block of $E(G)$ and since $B$ is reduced, and by \Cref{O2central} and (\ref{O2'}), we have $F^*(G)=Z(G)$, in which case $D \leq G=C_G(F^*(G)) \leq F^*(G)=Z(G)$, implying $D \lhd G$. Thus (\ref{E in out}) follows from (\ref{comptrans}) by the argument described in the proof of \cite[Lemma 5.6]{Elliot}. 
 
Following the proof of \cite[Lemma 5.6]{Elliot}, here considering $\Omega(Q)$ rather than all of $D$, this implies there is a subgroup isomorphic to $C_{2^n-1}$ in $C_{2^k-1} \wr C_t$ where $n=tk$. Considering exponents, \cite[Lemma 5.5]{Elliot} tells us there is no such embedding unless $t=1$ or $t=k$. If $t=k$, then we must that $\Omega(D) \cap L_i$ has rank $2$, with $R \leq \cap_{t}(L_i)$. This implies $R \leq Z(E(G)) \leq Z(G)$, and so the $B$-covered block of $E(G)$ would cover nilpotent blocks of every component, a contradiction. Thus $t=1$. 
\end{proof}

Before proceeding, we consider the classification of blocks of quasisimple groups with abelian defect groups, from \cite[Theorem 6.1]{ekks}. Here we present \cite[Corollary 4.8]{Elliot}, a table summarising that result and others from that paper for non-nilpotent blocks. 

\begin{theorem} \label{ekkssimple}
	Let $L$ be a quasisimple group, and $b$ be a non-nilpotent $2$-block of $\mathcal{O}L$ with abelian defect group $D$. Then the following table describes the six possible situations that can occur, including the isomorphism type of $L$, $D$, $\operatorname{Out}(L)$, and the inertial quotient $\mathbb{E}_b$ of $b$, where appropriate:

	\addtolength{\parskip}{-0.6cm}\begin{table}[H] \small
		\centering 
		\begin{tabular}{c c c c c l} 
			\hline\hline  
			\textbf{Group} 				& \textbf{Block}& $D$		& \footnotesize{Out($L)$}& $\mathbb{E}_b$ 	& Notes\\ [1ex]
			\hline 
			$SL_2(2^n)$ 	& Principal 	& $(C_2)^n$	& $C_n$		&	$C_{2^n-1}$	&	 $n \geq 2$	\\ [1ex]
			${}^2G_2(q)$& Principal 	& $(C_2)^3$	& $C_q$		&	$C_7 \rtimes C_3$	& $q=3^r$, $r$ odd  	\\ [1ex]
			$ J_1$ 		& Principal  	& $(C_2)^3$	& 	$1$	&	$C_7 \rtimes C_3$	&		\\ [1ex]
			$ Co_3$ 			& \small{Non-principal}	& $(C_2)^3$	& 	$1$	& 	$C_7 \rtimes C_3$	& \footnotesize{Unique non-principal}\\[1ex]
			\footnotesize{$L/Z(L)$ type $A_t$ or $E_6$}		& \footnotesize{Nilpotent-covered}& 	-	& 	-	&	-	&	\\ [1ex]
			$M_0 \times M_1 \leq L$ & $ b_0 \otimes b_1$	& 	$D_0 \times C_2^2$	&	-	& $C_3$ & \dagger \\ 		\hline 	\end{tabular}\end{table} \noindent \dagger { } In the final case, $L$ is of type $D_t(q)$ or $E_7(q)$ where $q$ is an odd prime power and $t/2$ is odd. Further, $b$ is Morita equivalent to $b_0 \otimes b_1$ where $b_i$ is a block of $M_i$; $M_0$ is abelian; $b_1$ is Puig equivalent to $\mathcal{O}A_4$ or $B_0(\mathcal{O}A_5)$; and $\mathbb{E}_b$ has order $3$ and centralises $D \cap M_0$. \addtolength{\parskip}{+0.6Cm}
\end{theorem}

\begin{remark} \label{allMor} \normalfont
 \cite[Theorem 1.5, Lemma 4.2 (xi)]{KoshitaniaBrouS} tells us that the non-principal block of $Co_3$ with defect group $(C_2)^3$ is Puig equivalent to the principal block of $\operatorname{Aut}(SL_2(8))$, and according to \cite[Example 3.3, Remark 3.4]{okuyamasomeexamples}, it had been essentially proved in \cite{Reee} that the principal block of $\operatorname{Aut}(SL_2(8))$ and the principal block of ${}^2G_2(q)$ are Puig equivalent (see \cite[Theorem 1.6]{KoshitaniaBrouS}).
\end{remark}

\section{Proof of the main theorem} \label{section proof}
Now we proceed towards proving \Cref{Main2} and \Cref{Maink}.

\begin{proposition}\label{O2central}
	Let $G$ be a finite group, and $B$ a $2$-block of $\mathcal{O}G$ with abelian defect group. Suppose that $O_2(G)$ is non-trivial and cyclic, and that $B$ covers a $G$-stable block $b_C$ of $C_G(O_2(G))$. Then $O_2(G)$ is central in $G$. 
\end{proposition}
\begin{proof}
	Since $D$ is abelian and $O_2(G) \leq D$, we have that $D \leq C_G(D) \leq C_G(O_2(G))$. Therefore, \Cref{alp15} (\ref{alp CG(D)<N}) tells us that ${b_C}^G=B$, and that $B$ is the unique block of $G$ covering $b_C$. Therefore, since $D \leq C_G(O_2(G))$, \Cref{alp15} (\ref{alpintersect}) tells us that $2 \not \divides [Stab_G(b_C):C_G(O_2(G))]$, and since $b_C$ is $G$-stable by assumption, this means $2 \not \divides [G:C_G(O_2(G))]$. 
	
	We know that $N_G(O_2(G))/C_G(O_2(G))$ embeds into $\aut(O_2(G))$. However, since $O_2(G)$ is a cyclic $2$-group, $\aut(O_2(G))$ is also a $2$-group. In particular, since $O_2(G) \lhd G$, this tells us that $G/C_G(O_2(G))$ is a $2$-group. Thus we must have that $[G:C_G(O_2(G))]=1$, i.e. that $G=C_G(O_2(G))$ and thus $O_2(G) \leq Z(G)$. 
\end{proof}

Recall that a group is \textit{supersolvable} if it has an invariant normal series where each factor is cyclic. 

\begin{lemma} \label{supersolvable}
	Let $G$ be a finite group, and $B$ a reduced $2$-block of $G$ with defect group $D$ and inertial quotient $\mathbb{E}$ such that $D=Q \times R$ with $Q \cong (C_{2^m})^n$ and $R \cong C_{2^r}$, and such that $\mathbb{E}$ centralises $R$ and contains an element of order $2^n-1$. If $D$ is not normal in $G$, then $G/F^*(G)$ is supersolvable and $G/L$ is solvable. 
\end{lemma}
\begin{proof}
	Applying \Cref{O2central} and 
	, we have that $F^*(G)=L Z(G)$. Since $G/F^*(G) \leq \out(F^*(G))$, this means we have $G/F^*(G) \leq \out(L)$. Further, $\out(L) = \out(L/Z(L))$, and Schreier's conjecture tells us this is solvable, implying $G/F^*(G)$ is solvable. Since, $F^*(G)/L$ is abelian and $G/F^*(G) \cong (G/L)/(F^*(G)/L)$, this tells us that $G/L$ is also solvable.	
	
	$B$ covers a block of a quasisimple group with defect group $D \cap L$, which means that we are in one of the six cases described by \Cref{ekkssimple}. Noting that $L$ cannot be type $D_4$, \cite[Table 5]{atlas} shows us that $\out(L)$ must be supersolvable. Since $G/F^*(G) \leq \out(L)$, we have that $G/F^*(G)$ is supersolvable. 
\end{proof}

Next we consider the index of $L$ in $G$. We recall a definition from \cite{higman1956length}: the lower $p$-series of a $p$-solvable group is $1 \leq O_{p'}(G) \leq O_{p',p}(G) \leq O_{p',p,p'}(G) \leq ...$ where $O_{p',p}(G)$ is the group containing $O_{p'}(G)$ such that $O_{p',p}(G)/O_{p'}(G)=O_p(G/O_{p'}(G))$, and so on. The $p$-length of $G$ is defined as the number of quotients in this series that are $p$-groups.


\begin{lemma} \label{2length}
	Let $G$ be a finite group, and $B$ a reduced block of $G$ with defect group $D$ and inertial quotient $\mathbb{E}$ such that $D=Q \times R$ with $Q \cong (C_{2^m})^n$ and $R \cong C_{2^r}$, and such that $\mathbb{E}$ centralises $R$ and contains an element of order $2^n-1$. If $D$ is not normal in $G$, then $Q \leq E(G)$. Further, either $R \leq E(G)$, or there exists a normal subgroup $N \lhd G$ containing $E(G)$, such that $2 \not \divides [N:E(G)]$ and $[G:N]=[R:R \cap E(G)] \divides 2^r$.
\end{lemma}
\begin{proof}
 Write $\overline{G}:=G/F^*(G)$. \Cref{supersolvable} tells us that $\overline{G}$ is supersolvable, so there is a chain $1=\overline{N}_0 \lhd \overline{N}_1 \lhd ... \lhd \overline{N}_{s-1} \lhd \overline{N}_s=G$ such that $[\overline{N}_i:\overline{N}_{i-1}]$ is a prime and $\overline{N}_i \lhd \overline{G}$, for each $i$. Taking the preimages of these subgroups under the canonical projection $G \to \overline{G}$, we get a corresponding chain $F^*(G)=N_0 \lhd N_1 \lhd ... \lhd N_s=G$ where $[N_i:N_{i-1}]$ is prime and $N_i \lhd G$, for each $i$. 
	
	Thus, for each $i$ we have that $[Q \cap N_i:Q \cap N_{i-1}]$ is $1$ or $2$. However, \Cref{Lemma} (\ref{DcapN}) tells us that $Q \cap N_i$ is trivial or of the form $(C_{2^{m'}})^{n}$ for some $m' \leq m$. Thus we must have $Q \cap N_i=Q \cap N_{i-1}$ for each $i$, implying that $Q \leq F^*(G)$. By \Cref{Lemma} (\ref{DcapN}), $Q \cap Z(G)=1$, and \Cref{O2central} tells us that $O_2(G) \leq Z(G)$, so we must have $Q \leq E(G)$. Note that \Cref{Lemma} (\ref{onecomp solv}) tells us that $E(G)=L$ is quasisimple.
	
	\Cref{solv quotient} tells us that $DL/L$ is a Sylow-$2$-subgroup of $G/L$, and thus if $[N_i:N_{i-1}]=2$ we must have $[R \cap N_i:R \cap N_{i-1}]=2$, implying that $RL/L$ is a Sylow-$2$-subgroup of $G/L$, and $[G:L]=[R:R \cap L] \leq 2^r$. Since $RL/L$ is abelian and $G/L$ is solvable, \cite[Chap $6$, ex $13$]{gorenstein} tells us the $2$-length of $G/L$ is $1$. That is, the lower $2$-series of $G/L$ is $$1 \leq O_{2'}(G) \leq O_{2',2}(G) \leq O_{2',2,2'}(G)=G/L$$  These are characteristic subgroups of $G/L$, and if we take preimages under the projection $G \to G/L$, we get a series $L \lhd N \lhd M \lhd G$ such that $L, M,$ and $N$ are characteristic subgroups of $G$, and $[M:N] \leq 2^r$, while $[G:M]$ \& $[N:L]$ are odd. Since $M/N \lhd G/N$, there's a non-trivial homomorphism $\frac{G/N}{M/N} \to  \out(M/N)$. However, $M/N$ is a cyclic $2$-group, so $\out(M/N)$ is a $2$-group, while $\frac{G/N}{M/N}  \cong G/M$ has odd order. This implies that $G=M$, so that we have a chain $L \lhd N \lhd G$ where $G/N \cong RL/L$, and $[N:L]$ is odd. \qedhere

\end{proof}


\begin{lemma} \label{main nonsplit}
Let $G$ be a finite group and $B$ a reduced block of $\mathcal{O}G$ with defect group $D=Q \times R$ and inertial quotient $\mathbb{E}$ such that $Q \cong (C_{2^m})^n$, $R \cong C_{2^r}$, and $\mathbb{E}$ centralises $R$ and contains an element of order $2^n-1$. Let $E(G) \lhd N \lhd G$ be the chain given in \Cref{2length}, and let $b_N$ be the block of $\mathcal{O}N$ covered by $B$. Suppose that $\mathbb{E}$ is cyclic. If $b_N$ is basic Morita equivalent to its Brauer correspondent in $N_N(D \cap N)$, then $B$ is Morita equivalent to its Brauer correspondent in $N_G(D)$. \end{lemma}

%

\begin{proof}
\Cref{Singers} (\ref{Singers autos}) tells us that $\mathbb{E}=E \rtimes F$ where $E \cong C_{2^n-1}$ and $F \leq C_n$. Since $\mathbb{E}$ is cyclic, we have $\mathbb{E}=E \cong C_{2^n-1}$. \Cref{2length} tells us that $Q \leq E(G)$. Thus, $b_N$ has defect group $D \cap N \cong Q \times R'$, where $R'\leq R \cong C_{2^r}$. Note that \Cref{IQ con} tells us $b_N$ has inertial quotient $E \cong C_{2^n-1}$, and so \Cref{Singers} (\ref{Singers trans}) tells us that $E$ acts freely on the non-trivial elements of $[D,E]=Q$. 
	
First, if $R \leq E(G)$, then we have $G=N$ and $B=b_N$, and so we are done by assumption.
	
	Next, suppose that $1 < R \cap E(G) < R$, so that $D$ is a non-split extension of $D \cap N$. Let $B_D$ be the Brauer correspondent of $B$ in $N_G(D)$. To show that $B$ and $B_D$ are basic Morita equivalent, we will show there is a Morita equivalence between $B$ and $B_D$ that is compatible with the $*$-structure. To do so, it is sufficient by \cite[Propositions 2.7 \& 3.6]{wuzzz} to show that there is a Morita equivalence between $B$ and $\mathcal{O}(D \rtimes E)$ that is compatible with the $*$-structure.
	
	When $r=0$, we have $D=Q$ and $G=N$, so $B$ is basic Morita equivalent to $B_D$ by assumption. By \cite[Proposition 3.6]{wuzzz}, this implies there is a Morita equivalence between $B$ and $B_D$ that is compatible with a local system, and by \cite[Proposition 2.7]{wuzzz}, this is equivalent to the existence of a Morita equivalence between $B$ and $\mathcal{O}(Q \rtimes E)$ that is compatible with the $*$-structure. 
	
	We proceed by induction on $r$. Since $b_N$ has defect group $D \cap N=Q \times C_{2^{r'}}, 1 \leq r' < r$ and inertial quotient $E$, by induction there is a Morita equivalence between $\mathcal{O}Nb_N$ and $\mathcal{O}((D \cap N) \rtimes E)$ that is compatible with the $*$-structure as in \cite{star}. Since $\mathbb{E}$ is cyclic by assumption, and we know $C_Q(E)=1$, \cite[Proposition 5.10]{wuzzz} then tells us there is a Morita equivalence compatible with the $*$-structure between $\mathcal{O}GB$ and $\mathcal{O}(D \rtimes E)$, as claimed. 
	
	Finally, suppose that $R \cap E(G)=1$, so that $b_N$ is basic Morita equivalent to $\mathcal{O}(Q \rtimes E)$. Note that $\mathcal{O}(Q \rtimes E)$ has $|E|$ simple modules, and since basic Morita equivalence preserves the number of simple modules, we have $l(b_N)=|E|$. Further, as $D$ is abelian and $b_N$ is $G$-stable, \cite[Proposition 3.1]{eatonlivesyDons} tells us that $D$ acts as inner automorphisms on $b_N$. Since $G=N \rtimes R$, this implies that $l(B)=l(b_N)=|E|$. Thus, writing $\mathcal{N}:=N_G(D,B_D)$, \Cref{watanabepis} tells us there exists a perfect isometry $\kappa: \PI(\mathcal{N},{B_D}^\mathcal{N}) \to \PI(G,B)$ such that $\kappa(\lambda * \zeta) = \lambda * \kappa(\zeta)$ for $\lambda \in \operatorname{Irr(C_D(\mathcal{N}))}$ and $\zeta \in \PI(\mathcal{N},{B_D}^\mathcal{N})$. 
	
	
	Now, ${B_D}^\mathcal{N}$ has inertial quotient $E$, and thus we can apply the argument above again to get a perfect isometry $\epsilon:  \PI(\mathcal{N},{B_D}^{\mathcal{N}}) \to \PI(C_G(R),B_R)$ such that $\epsilon(\lambda * \zeta)$ for all $\lambda \in \operatorname{Irr}(C_D(\mathcal{N}))$ and $\zeta \in \PI(\mathcal{N},{B_D}^\mathcal{N})$.  Thus, the composition $\kappa \circ \epsilon^{-1}$ satisfies the conditions of \Cref{watanabesplitting}, and so $B \cong b_N \otimes_\mathcal{O} \mathcal{O}R$ as $\mathcal{O}$-algebras. Thus $B$ is Morita equivalent to $\mathcal{O}((Q \rtimes E) \times R)$. \end{proof} 



\begin{lemma} \label{main kosh}
Let $G$ be a finite group and $B$ a block of $kG$ with defect group $D=Q \times R$ and inertial quotient $\mathbb{E}$ such that $Q \cong (C_{2^m})^n$, $R \cong C_{2^r}$, and $\mathbb{E}$ centralises $R$ and contains an element of order $2^n-1$. Suppose that $\mathbb{E}$ is non-cyclic and that $R \cap E(G) =1$. If $B$ covers an inertial block of $E(G)$, then $B$ is Morita to its Brauer correspondent in $kN_G(D)$. 
\end{lemma}
\begin{proof}
\Cref{Lemma} (\ref{onecomp solv}) tells us that $E(G)=L$ is quasisimple, and \Cref{2length} tells us there is a chain $L \lhd N \lhd G$ such that $[N:L]$ is odd and, since $R \cap L=1$, that  $G/N \cong R$. In particular, $D=(D \cap N) \times R$, and so \Cref{koshkul} tells us that $B \cong b_N \otimes_k kR$ as $k$-algebras. Thus $B$ is  Morita equivalent to $k(Q \rtimes \mathbb{E}) \otimes_k kR$, and thus to $k((Q \rtimes \mathbb{E}) \times R)$. 
\end{proof}

Now we prove \Cref{Main2}. 

%
\begin{proof}
Let $G$ be a finite group and $B$ a block of $\mathcal{O}G$ with defect group $D=Q \times R$ and inertial quotient $\mathbb{E}$ such that $Q \cong (C_{2^m})^n$, $R \cong C_{2^r}$, and $\mathbb{E}$ centralises $R$ and contains an element of order $2^n-1$, with $n,m,r \geq 1$.
	
\Cref{Singers} (\ref{Singers autos}) tells us that $\mathbb{E}=E \rtimes F$, where $E \cong C_{2^n-1}$ and $F$ is isomorphic to an odd order subgroup of $C_n$, such that $E$ acts freely on $Q \backslash \{1\}$ and transitively on $\Omega(Q) \backslash \{1\}$, and a generator of $F$ acts on $\Omega(Q)$ as the $p^{s}$th-powering, where $|F|=n/s$. 

We suppose $B$ is a minimal counterexample: in particular, suppose $([G:O_{2'}(Z(G))],|G|)$ is minimised in the lexicographic ordering such that $B$ does not lie in any case of \Cref{Main2}. We may assume $B$ is reduced in the manner described in \Cref{reduce remark}: suppose $N \lhd G$ and $b_N$ is a block of $N$ covered by $B$. \Cref{F1} tells us that $B$ covers a block $B_I$ of $\stab_G(b_N)$ that is basic Morita equivalent to $B$. By minimality, we must have $B=B_I$, so $b_N$ is $G$-stable. Applying this argument to each normal subgroup, $B$ is quasiprimitive. Further, by minimality and \Cref{F2 cor}, if $N \lhd G$ and $B$ covers a nilpotent block of $N$, then $N \leq O_2(G)Z(G)$.

Note that the Schur multiplier of $\mathbb{E}$ is trivial by \cite[Proposition 5.3]{Elliot}. Thus, if $D$ is normal in $G$, then \Cref{kuls} tells us that $B$ is Puig equivalent to $\mathcal{O}(D \rtimes \mathbb{E})$, placing us in case (\ref{2Inert}) and contradicting our assumption that $B$ is a counterexample. Therefore $D$ is not normal in $G$, and so \Cref{Lemma} (\ref{O2}) tells us that $O_2(G) \leq R$, and \Cref{Lemma} (\ref{O2'}) implies that $O_{2'}(G) \leq Z(G) \leq Z(F^*(G)) \leq O_{2}(G)Z(G)$. \Cref{O2central} tells us that $O_2(G) \leq Z(G)$, and so, by \Cref{Lemma} (\ref{onecomp solv}), $G$ contains a single component $L \lhd G$, $F^*(G)=LZ(G)$, $G/F^*(G) \leq \out(L)$, and \Cref{supersolvable} tells us that $G/F^*(G)$ is supersolvable. Finally, \Cref{2length} tells us that $Q \leq L$, and that there is a chain $L \lhd N \lhd G$ such that $2 \not \divides [N:L]$, and $[G:N]=[R:R \cap L] \divides 2^r$. Let $b$, $b^*$ and $b_N$ denote the $B$-covered blocks of $L$, $F^*(G)$, and $N$ respectively, and note that $b^*$ covers $b$ and that $b$ and $b_N$ share a defect group. 

Consider $\overline{G}:= G/O_{2'}(G)$, and note that there is a block $\overline{B}$ of $\overline{G}$ whose defect group and inertial quotient are isomorphic to those of $B$. However, $[G:O_{2'}(G)]=[\overline{G}:O_{2'}(\overline{G})]$, while $|G|>|\overline{G}|$. Thus, by minimality, we have that $G=\overline{G}$, and so $F^*(G)=LZ(G)=LO_2(G)$.

Next, we consider which candidate structures for $b$ and $L$ from \Cref{ekkssimple} are compatible with what we now know. Let $\eta$ denote the rank of $D \cap L$, so that $\eta=n$ if $R \cap L=1$, and $\eta =n +1$ otherwise.  \begin{itemize}
	
\item Suppose $L \cong SL_2(2^{\eta})$, in which case $D \cap L \cong (C_2)^{\eta}$, and $b$ has inertial quotient $E_b \cong C_{2^{\eta}-1}$. We claim that $\eta=n$. Since $G/F^*(G) \leq \out(SL_2(2^\eta)) \cong C_\eta$, and $F^*(G)/L \leq Z(G)$, \cite[Theorem 4.12]{Elliot} tells us that there is some subgroup of $E_b$ that is isomorphic to a normal subgroup of $\mathbb{E}$. However, $E_b$ is generated by a Singer cycle in $GL_\eta(2)$, and thus no non-trivial elements of $E_b$ fix any non-trivial elements of $D \cap L$ by \Cref{Singers} (\ref{Singers trans}). If $R \cap L \neq 1$, this implies there is a subgroup of $\mathbb{E}$ that does not centralise $R$. This is a contradiction, and thus $\eta=n$. 

Next, we claim that $F^*(G)=L \times R$. Consider the action of $RF^*(G)/F^*(G)$ on $D \cap F^*(G)$. We have that $RF^*(G)/F^*(G) \leq G/F^*(G) \leq \out(SL_2(2^n))$, and the outer automorphisms of $SL_2(2^n)$ are the field automorphisms of $SL_2(2^n)$. In particular, as in \cite[Lemma 4.10]{Elliot}, there is a defect group $S$ of $b$ such that $\out(SL_2(2^n))$ centralises no non-trivial elements of $S$. However, $D$ is abelian, so $RF^*(G)/F^*(G)$ must centralise the defect groups of $b$. This implies $R \leq F^*(G)$, and so $R=Z(G)$ and $F^*(G)=SL_2(2^n) \times R$. 

Since $G=ND$, this implies that $G=N \times R$. Thus, we have $F^*(G)=F^*(N) \times R$, and so $N \leq \aut(L)$. Since \Cref{IQ con} tells us that $b_N$ has inertial quotient $\mathbb{E}$, \cite[Lemma 4.10]{Elliot} tells us $[N:L]=[\mathbb{E}:E_b]$. That is, $G \cong (SL_2(2^n) \rtimes F) \times R$ where $\mathbb{E}=E_b \rtimes F$. This places us in case (\ref{2SL}), contradicting our assumption that $B$ is a counterexample.

\item Suppose $L \cong {}^2G_2(q)$, where $q=3^f$ for $f$ odd. Note that $\out({}^2G_2(q)) \cong C_f$, so $[G:F^*(G)]$ is odd. Thus $R \leq F^*(G)=L \times Z(G)$, and so either $R \leq L$ or $R=Z(G)$. If $R \leq L$, then $G \leq \aut(L)$ and \cite[Proposition 3.1]{eaton3} tells us that $B$ is Puig equivalent to $b$. This implies $B$ has inertial quotient $C_7 \rtimes C_3$, which does not centralise $R$ - a contradiction. Therefore, $F^*(G)=L \times R$. This implies $F^*(G)= F^*(N) \times R$ and so $N \leq \aut(L)$. Thus, $G=N \times R$, where $N \leq \aut({}^2G_2(q))$. Note that \cite[Proposition 3.1]{eaton3} tells us that $b_N$ is Puig equivalent to $b$, implying by \Cref{IQ con} that $\mathbb{E} \cong C_7 \rtimes C_3$. Further, by \Cref{allMor}, this implies $b_N$ is Puig equivalent to the principal block of $\aut(SL_2(8))$, and since $G$ is a direct product, $B$ is M we have that $B \cong b_N \otimes_\mathcal{O} \mathcal{O}R$, placing us in case (\ref{2SL}). 

\item Suppose $L \cong Co_3$. Since $\operatorname{Out}(Co_3)=1$, we have $G=F^*(G) \cong Co_3 \times O_2(G)$. We must have $R \cap Co_3=1$, since otherwise the inertial quotient of $B$ does not centralise $R$. Thus, $G \cong Co_3 \times R$, and so by \Cref{allMor} and as above, $B$ is Puig equivalent to the principal block of $\aut(SL_2(8)) \times R$, placing us in case (\ref{2SL}) and giving a contradiction. 
	
\item Suppose $L \cong J_1$. Since $\operatorname{Out}(J_1)=1$, we have $G=F^*(G) \cong J_1 \times O_2(G)$. We must have $R \cap J_1=1$, since otherwise the inertial quotient of $B$ does not centralise $R$. Thus, $G \cong J_1 \times R$, placing us in case (\ref{2J1}) and providing a contradiction.
	
	\item Suppose $b$ is nilpotent-covered. Then \Cref{zhou} (\ref{zhou b nilp b inert}) tells us $b$ is inertial, and \Cref{zhou} (\ref{zhou b inert B inert}) then tells us $b_N$ is inertial. If $R \leq L$ then $G=N$, and so $B$ is inertial. If $\mathbb{E}$ is cyclic, then \Cref{main nonsplit} tells us that $B$ is inertial. Thus if $R \leq L$ or $\mathbb{E}$ is cyclic, we are in case (\ref{2Inert}), providing a contradiction.

\item Suppose $b$ lies in the final case of \Cref{ekkssimple}, so there is a decomposition $D \cap L=(D \cap M_0) \times (D \cap M_1)$ where $D \cap M_1 \cong (C_{2})^2$. Clearly $Q \cap M_1 \neq 1$, and by considering the rank of $D$ and $D \cap L$, this implies $Q$ is elementary abelian. \cite[Propositions 5.3]{ekks} tells us that $L$ is of type $D_t(q)$ or $E_7(q)$ with $t/2$ odd and $q$ a power of an odd prime. The $2$-length argument in \Cref{2length} applies  to $G/F^*(G)$, and so there exists a chain $F^*(G) \lhd M \lhd G$ such that $M/F^*(G)$ is odd, and $[G:M]=[D:D \cap F^*(G)]$. Let $b_M$ denote the block of $M$ covered by $B$, and note that \Cref{IQ con} tells us that $b_M$ has inertial quotient $\mathbb{E}$. Note that $M/F^*(G)$ is an odd subgroup of $\out(L)$, and by examining \cite[Table 5]{atlas}, this implies that $M/F^*(G)$ is cyclic. Thus, since $F^*(G)/L$ is abelian, \cite[Theorem 4.12]{Elliot} tells us that $E_b \lhd \mathbb{E}$. However, $E_b$ centralises all $D \cap M_0$, and so if $n>2$ we have a contradiction to \Cref{Singers} (\ref{Singers auto nofix}).

Thus we have $n=2$, and \cite[Theorem 3.2]{eatonlivesrank3} tells us that $B$ is Morita equivalent to $\mathcal{O}(A_4 \times R)$ or $B_0(\mathcal{O}(A_5 \times R))$. The former places us in case (\ref{2Inert}), while the latter places us in case (\ref{2SL}), each giving a contradiction. \end{itemize}

Since all the cases above lead to a contradiction, $B$ cannot be a minimal counterexample, and so \Cref{Main2} must hold. \qedhere\end{proof} We conclude by proving \Cref{Maink}.

%

\begin{proof}
Let $G$ be a finite group and $\tilde{B}$ a block of $kG$ with defect group $D=Q \times R$ and inertial quotient $\mathbb{E}$ such that $Q \cong (C_{2^m})^n$, $R \cong C_{2^r}$, and $\mathbb{E}$ centralises $R$ and contains an element of order $2^n-1$, with $n,m,r \geq 1$. The arguments in the proof of \Cref{Main2} all apply to $\tilde{B}$ as they did to the block of $\mathcal{O}G$, with the sole exception of applying \Cref{main nonsplit}. To replace this, suppose that either $R \cap E(G)=1$ or $R \leq E(G)$, and that $\tilde{B}$ covers an inertial block of the normal subgroup $N$ defined in \Cref{2length}. If $R \leq L$ then $N=G$ and so $\tilde{B}$ is inertial, while if $R \cap L=1$ then \Cref{main kosh} tells us that $\tilde{B}$ is Morita equivalent to $k((Q \rtimes \mathbb{E}) \times R)$. 
\end{proof}

\section*{Acknowledgements}
This paper is part of the work towards my PhD at the University of Manchester, which is supported by a University of Manchester Research Scholar award. I am deeply grateful to Charles Eaton, my PhD supervisor, for his crucial advice and patient support. I'm thankful to Michaels Livesy for his useful input, and to Cesare Ardito for countless productive discussions.

\small
\bibliographystyle{amsplainab}\bibliography{biblio}

\end{document}